\documentclass[12pt,jmp]{revtex4-1}
\usepackage[margin=1in]{geometry}  
\usepackage{graphicx}              
\usepackage{amsmath}               
\usepackage{amsfonts}              
\usepackage{amsthm}                

\usepackage{color}

\def\Re{\mathbb{R}}
\newtheorem{claim}{Claim}[section]
\newtheorem{thm}{Theorem}[section]

\newtheorem{lem}{Lemma}[section]

\newtheorem{cor}{Corollary}[section]
\newcommand{\fim}{\hfill\rule{2mm}{2mm}}

\numberwithin{equation}{section}

\usepackage{graphicx}
\usepackage{dcolumn}
\usepackage{bm}

\begin{document}

\preprint{AIP/123-QED}

\title[A critical  nonlinear fractional  elliptic equation  ]{A critical  nonlinear fractional  elliptic equation with saddle-like potentical in $\mathbb{R}^N$ }
\thanks{Research of C. O. Alves partially supported by  CNPq/Brazil  304036/2013-7  and INCTMAT/CNPq/Brazil. \ O.H.Miyagaki is corresponding author, and He  was partially supported by INCTMAT/CNPq/Brazil, CNPq/Brazil 304015/2014-8  and CAPES/Brazil 
Proc. 2531/14-3.  }

\author{Claudianor  O. Alves}
 \affiliation{Universidade Federal de Campina Grande, \\  Unidade Acad\^emica de Matem\'atica,\\  CEP: 58429-900 - Campina Grande-PB, Brazil }
 \email{ coalves@dme.ufcg.edu.br}

\author{Olimpio H. Miyagaki}
\affiliation{Departmento de Matem\'atica, \\
 Universidade Federal de  Juiz de Fora,\\ 
 CEP:  36036-330 - Juiz de Fora-MG, Brazil
}
\email{ohmiyagaki@gmail.com}

\date{\today}

\begin{abstract}
In this paper, we study  the existence of positive solution for the following class of fractional elliptic equation 
$$
\epsilon^{2s} (-\Delta)^{s}{u}+V(z)u=\lambda |u|^{q-2}u+|u|^{2^{*}_{s}-2}u\,\,\, \mbox{in} \,\,\, \mathbb{R}^{N},
$$
where  $\epsilon, \lambda >0$ are positive parameters, $q \in (2,2^{*}_{s}), 2^{*}_{s}=\frac{2N}{N-2s}, $  $N > 2s,$ $s \in (0,1),$  $ (-\Delta)^{s}u$ is the fractional laplacian, and $V$ is a saddle-like potential. The result is proved by using minimizing method constrained to the Nehari manifold. A special minimax level is obtained by using an argument made  by Benci and Cerami.
\end{abstract}



\maketitle

\vspace{0.5 cm}


\section{Introduction}

The nonlinear fractional  Schr\"{o}dinger
equation

$$
i\epsilon \displaystyle \frac{\partial \Psi}{\partial t}=\epsilon^{2s}(-\Delta)^s
\Psi+(V(z)+E)\Psi-f(\Psi)\,\,\, \mbox{for all}\,\,\, z \in
\Re^N,\eqno{(NLS)}
$$
where $N > 2s,$ $\epsilon > 0$, $V,f$ are continuous functions, has been studied in recent years by many researchers. The standing waves solutions of $(NLS)$, namely, $\Phi(z,t)=exp(-i Et)u(z),$ where $u$ is a solution of the
fractional elliptic equation
$$
\ \  \left\{
\begin{array}{l}
\epsilon^{2s} (- \Delta)^su + V(z)u=f(u)
\ \ \mbox{in} \ \ \Re^N,
\\
u \in H^s(\Re^N),\quad u > 0 \ \ \mbox{on} \ \ \Re^N.
\end{array}
\right.
\eqno{(P_{\epsilon})}
$$
In the local case, that is, when $s=1,$ the existence and concentration of
positive solutions for general semilinear elliptic equations
$(P_{\epsilon})$, for the case $N \geq 2$,  have been extensively
studied,  for example,  by  \cite{AOS,Alves10,BPW,BW0,Pino,PFM,OS,FW,Oh1,R,W},  and their references. Rabinowtz in \cite{R},  proved the
existence of positive solutions of $(P_{\epsilon})$, for $\epsilon > 0$
small, imposing a global condition,
\[
\liminf_{|z| \rightarrow \infty} V(z) > \inf_{z \in
	\mathbb{R}^N}V(z)=\gamma >0.
\]
 In fact, these solutions concentrate at 
global minimum points of $V$ as  $\epsilon$ tends to 0, c.f.  Wang in \cite{W}.  del Pino and Felmer in  \cite{Pino}, assuming a local condition, namely,  there is an open and bounded set $\Lambda$ compactly contained  in $\mathbb{R}^{N}$ satisfying
$$
0< \gamma \leq V_0 =\inf_{z\in\Lambda}V(z)< \min_{z \in
	\partial\Lambda}V(z), \eqno{(V_{1})}
$$
 established  the existence of positive  solutions which concentrate around local minimum of $V,$  by introducing a penalization method. 

In \cite{PFM}, del Pino, Felmer and Miyagaki considered the case where potential $V$ has a   saddle like geometry. They assumed that $V$ is bounded and  $V \in C^{2}(\mathbb{R}^{N})$, verifying the following conditions: Fixed two subspaces $X,Y  \subset \mathbb{R}^{N}$ such that
$
\mathbb{R}^{N}=X \oplus Y,
$  
define  $c_0,c_1>0$  given by \\
$$
\displaystyle c_0=\inf_{z \in \mathbb{R}^{N}}V(z)>0
\quad \mbox{
and}\quad 
c_1=\displaystyle \sup_{x \in X}V(x),
$$
satisfying  the following geometric condition  \\

\noindent $(V_1)$  
$$
c_0=\inf_{R>0}\sup_{x \in \partial B_R(0) \cap X }V(x)<\inf_{y \in Y}V(y).
$$
In addition to the above hypotheses, they imposed the conditions below:\\

\noindent $(V_2)$ \quad The functions $V, \frac{\partial V}{\partial x_i}$ and $\frac{\partial^{2} V}{\partial x_i \partial x_j}$ are bounded 
in $\mathbb{R}^{N}$ for all $i,j \in \{1,...,N\}$.  \\

\noindent $(V_3)$ \quad $V$ satisfies the Palais-Smale condition, that is, if $(x_n) \subset \mathbb{R}^{N}$ is a sequence such that $(V(x_n))$ is bounded and $\nabla V(x_n) \to 0$, then $(x_n)$ possesses a convergent subsequence in $\mathbb{R}^{N}$. \\

Using the above conditions on $V$, and supposing that 
$$
c_1<2^{\frac{2(p-1)}{N+2-p(N-2)}}c_0,
$$
the authors in \cite{PFM} showed  the existence of positive solutions for  $(P_{\epsilon})$ with $ f(u)=|u|^{p-2}u$,
where $p \in (2,2^{*})$ if $N \geq 3$ and $p \in (2,+\infty)$ if $N=1,2$, for $\epsilon>0$ small enough. Here  $2^*=\frac{2N}{N-2} $ is the critical Sobolev exponent.

Motivated by the results obtained in \cite{PFM}, with the potential $V$ having the same geometry as  considered in \cite{PFM}, Alves in \cite{AlvesNovo}  proved  the existence of positive solution result for  $(P_{\epsilon})$, not only  with  $f$ having  an exponential critical growth, for $N=2,$  but  also  with $f(u)=|u|^{q-2}u+|u|^{2^{*}-2}u $, \  where  $q \in (2,2^{*})$ and    $N\geq 3.$

\vspace{0.3 cm}
In the nonlocal case, that is, when $s\in (0,1),$ even in the subcritical case,  there are only few references on existence and/or concentration phenomena for fractional nonlinear equation $(P_{\epsilon})$, maybe because techniques developed for local case can not be adapted imediately, c.f. \cite{Secchi}. We would like to cite \cite{Secchi, Shang} for the  existence of positive solution, imposing a global condition on $V.$   In \cite{Davila} is studied the existence and concentration phenomena for potential verifying local condition $(V_1)$, and in \cite{Chen, ShangZhang} a concentration phenomenon is treated  near of non degenerate critical point of $V.$

\vspace{0.3cm}
 By using  the same approach as  in \cite{AlvesNovo}, we will establish an existence result of positive solution for the following class of problem involving critical Sobolev  exponents
$$
\epsilon^{2}(-\Delta)^{s}u+V(x)u=\lambda |u|^{q-2}u+|u|^{2^{*}_s-2}u\,\,\, \mbox{in} \,\,\, \mathbb{R}^{N}, \eqno{(P_\epsilon)_*}
$$
where $s\in (0,1),$  $\epsilon, \lambda >0$ are positive parameters, $q \in (2,2_s^{*}), 2_s^{*}=\frac{2N}{N-2s}, N >2s,$ with  $V$   verifying the  above conditions and a relation on the numbers $c_0$ and $c_1,$ given by

\vspace{0.3 cm}

\noindent $(V_4)$ \quad $m_\lambda(V(0)) \geq 2m_\lambda(c_0)$ and $c_1 \leq c_0 + \frac{3}{5}\left(\frac{1}{2}-\frac{1}{q}\right)c_0,$

\vspace{0.5 cm}

\noindent where  $m_\lambda(A)$ is the   mountain pass level of the functional

$$
J_{\lambda, A}(u)=\frac{1}{2}\int_{\mathbb{R}^{N}}|\xi|^{2s}|\widehat{u}|^2 d\xi+  \frac{A}{2}\int_{\mathbb{R}^{N}}|u|^{2} \,dx- \frac{\lambda}{q} \int_{\mathbb{R}^{N}} |u|^{q} dx-\frac{1}{2^{*}_{s}}\int_{\mathbb{R}^{N}}|u|^{2^{*}_{s}} dx
$$
 defined in $H^{s}(\mathbb{R}^{N})$.  It is well known that the equality below   holds
$$
m_\lambda(A)=\inf_{u \in H^{s}(\mathbb{R}^{N}) \setminus \{0\}}\left\{\max_{t \geq 0}J_{\lambda,A}(tu) \right\}.
$$

Using the above notation we are able to state our main result 

\begin{thm} \label{T2} Assume that $(V_1)-(V_4)$ hold. Then,  there is $\epsilon_0>0$,  and $\lambda^{*}>0$ independent of $\epsilon_0$, such that   $(P_\epsilon)_*$ has a positive solution for all $\epsilon \in (0,\epsilon_0]$ and  $\lambda \geq \lambda^{*}$. 
\end{thm}

The proof of  Theorem \ref{T2} was inspired from    \cite{AlvesNovo} and  \cite{PFM},  however  since we are working 
with fractional laplacian,  the estimates for this class of  nonlocal problem are very delicate and different from  those used in the  local problems. We minimize the energy function constrained on the  Nehari manifold, and to this end, we modify the barycenter properly for our problem

We recall that, for any $s \in (0,1)$, the fractional Sobolev space $H^{s}(\Re^N)$ is definied by
\[
H^{s}(\Re^N)=\Big\{u\in L^2(\Re^N): \ \int_{\Re^{2N}}\frac{(u(x)-u(y))^2}{|x-y|^{N+2s}}\ d x\ d y<\infty\Big\},
\]
endowed with the norm
$$
\|u\|_{H^s(\Re^N)}=\Big(|u|_{L^2(\mathbb{R}^{N})}^2+\int_{\Re^{2N}}\frac{(u(x)-u(y))^2}{|x-y|^{N+2s}}\ d x\ d y\Big)^{1/2}.
$$
The fractional  Laplacian, $(-\Delta)^{s}u,$  of a smooth function $u:\Re^{N} \rightarrow  \Re$ is defined  by  
$$
{\mathcal F}((-\Delta)^{s}u)(\xi)=|\xi|^{2s}{\mathcal F}(u)(\xi), \ \xi \in \Re^N,
$$ 
where ${\mathcal F}$ denotes the Fourier transform, that is, 
\[
{\mathcal F}(\phi)(\xi)=\frac{1}{(2\pi)^{\frac{N}{2}}} \int_{\mathbb{R}^N} \mathit{e}^{-i \xi \cdot x} \phi (x)  \, \ d x \equiv \widehat{\phi}(\xi) ,  
\]
for functions $\phi$ in  the  Schwartz class.
Also $(-\Delta)^{s}u$  can be equivalently represented \cite[Lemma 3.2]{nezza} as
$$
(-\Delta)^{s} u(x) = -\frac{1}{2} C(N,s)\int_{\Re^N}\frac{(u(x+y)+u(x-y)-2 u(x))}{|y|^{N+2s}}\ d y, \ \forall x \in \Re^N,
$$
where 
$$C(N,s)=(\int_{\Re^N}\frac{(1- cos\xi_1)}{|\xi|^{N+2s}}d\xi)^{-1},\  \xi=(\xi_1,\xi_2,\ldots,\xi_N).$$
Also, in light of \cite[Propostion~3.4,Propostion~3.6]{nezza}, we have
\begin{equation}
\label{equinorm}
|(-\Delta)^{s/2} u|^2_{L^2(\Re^N)}=\int_{\Re^N}|\xi|^{2s}|\widehat{u}|^2d\xi=\frac{1}{2}C(N,s)\int_{\Re^{2N}}\frac{(u(x)-u(y))^2}{|x-y|^{N+2s}}\  d x \ d y, \quad \text{for all $u\in H^{s}(\Re^N)$},
\end{equation}
and, sometimes, we identify these two quantities by omitting the normalization constant $\frac{1}{2} C(N,s).$
For $ N > 2s,$  from \cite[Theorem 6.5]{nezza} we also know that, for any $p \in [ 2, 2^{*}_{s}]$,
there exists $C_p>0$ such that
\begin{equation}
\label{emb}
|u|_{L^p(\mathbb{R}^{N})}\leq C_p\|u\|_{H^{s}(\Re^N)},
\,\quad\text{for all $u\in H^{s}(\Re^N)$}.
\end{equation}

The best Sobolev constant $S$ is given by (see \cite{Coti})
$$S=\inf_{u \in H^{s}_{0}(\Re^N)\setminus \{0 \}}\frac{\int_{\Re^{2N}} \frac{|u(x)-u(y)|^2}{|x-y|^{N+2s}}dx dy}{(\int_{\Re^N}|u|^{2^{*}_{s}}dx)^{\frac{2}{2^{*}_{s}}  }  },$$
which is attained by
$$v_0(x)=\frac{c}{(\theta^2 +|x-x_0|^2)^{\frac{N-2s}{2}}}, \ x \in \Re^N, $$
where $c, \theta >0$ are constants and $x_0\in \Re^N$ fixed, and 
$$H^s_{0}(\mathbb{R}^{N})=\{u \in L^{2^{*}_{s}}(\Re^N): \ |\xi|^s\widehat{u} \in L^2(\Re^N)\}.$$

Before to conclude this introduction, we would like point out that using the change variable $v(x)=u(\epsilon x)$, it is possible to prove that $(P_\epsilon)$ is equivalent to the following problem
$$
(-\Delta)^s {u}+V(\epsilon x)u=f(u) \,\,\, \mbox{in} \,\,\, \mathbb{R}^{N}, \eqno{(P_\epsilon)'}
$$
where 
$$
f(t)=\lambda |t|^{q-2}t+|t|^{2^{*}_{s}-2}t, \quad \forall t \in \Re.
$$

In the present paper, we denote by $I_\epsilon$ the energy functional associated with $(P_\epsilon)'$ given by
\begin{eqnarray*}
I_\epsilon(u)&=&\frac{1}{2}\int_{\mathbb{R}^{N}}|\xi|^{2s}|\widehat{ u}|^{2}d\xi+\frac{1}{2}\int_{\mathbb{R}^{N}}V(\epsilon x)|u|^{2}\,dx- \frac{\lambda}{q} \int_{\mathbb{R}^{N}} |u|^{q} dx- \frac{1}{2^{*}_{s}}\int_{\mathbb{R}^{N}}|u|^{2^{*}_{s}} dx, 
\end{eqnarray*}
 for all $ u \in H^{s}(\mathbb{R}^{N}).$
It is standard to prove that $I_\epsilon \in C^1(H^{s}(\Re^N),\Re)$  with Gateaux derivative 
\begin{eqnarray*}
I'_\epsilon(u)v&=&\int_{\mathbb{R}^{N}}|\xi|^{2s}\widehat{ u}\  \widehat{ v}d\xi+\int_{\mathbb{R}^{N}}V(\epsilon x)u v\,dx- \lambda \int_{\mathbb{R}^{N}} |u|^{q-2}u v dx- \int_{\mathbb{R}^{N}}|u|^{2^{*}_{s}-2} uvdx
\end{eqnarray*}
 for all $ u, v  \in H^{s}(\mathbb{R}^{N}).$
This way, $u \in H^{s}(\mathbb{R}^{N})$ is a weak solution for $(P_\epsilon)'$ if, and only if, $u$ is a critical point for $I_\epsilon$.  

\vspace{0.5 cm}

\noindent \textbf{Notation:} In this paper we use the following
notations:
\begin{itemize}
	\item  The usual norms in $L^{t}(\mathbb{R}^{N})$ and $H^{s}(\mathbb{R}^{N})$ will be denoted by
	$|\,. \,|_{t}$ and $\|\;\;\;\|$ respectively.

	\item   $C$ denotes (possible different) any positive constant.
	
	\item   $B_{R}(z)$ denotes the open ball with center at $z$ and
	radius $R$.

\end{itemize}

\section{Technical results}

The next lemma is a Lions Lemma type result  which
can be adapted to our case, see \cite[Proposition II.3]{Secchi}.

\begin{lem} \label{lions} Let $(u_n) \subset H^{s}(\mathbb{R}^N)$ be a sequence  such that,

if there is $R>0$ such that
$$
\lim_{n \to +\infty}\sup_{z \in \mathbb{R}^N}\int_{B_R(z)}|u_n|^{2}\,dx=0,
$$
then
$$
\lim_{n \to +\infty}\int_{\mathbb{R}^N}|u_n|^q\,dx=0,\ \quad \forall q \in (2, 2^{*}_{s}).
$$
\end{lem}

As a consequence of the above lemma, we have the following result
\begin{cor} \label{sequencia}  Let $(u_n) \subset H^{s}(\mathbb{R}^N)$ be a  $(PS)_c$ sequence
 for $I_\epsilon,$  that is, 
$$ I_{\epsilon}( u_n) \to c,\quad \mbox{and}\quad I^{´}_{\epsilon}(u_n) \to  0,$$
with $0<c< \frac{s}{N} (2^{-1}C(N,s)S)^{N/2s}$ and $u_n \rightharpoonup 0.$
Then, there exists $(z_n) \subset \mathbb{R}^N$ with $|z_n| \to +\infty$ such that
$$
v_n=u_n(\cdot+z_n) \rightharpoonup v \not=0 \quad \mbox{in} \quad H^{s}(\mathbb{R}^N).
$$
\end{cor}
\noindent {\bf Proof.}  We claim  that for any $R>0$, 
$$
\lim_{n \to +\infty}\sup_{z \in \mathbb{R}^N}\int_{B_R(z)}|u_n|^{2}\,dx>0.
$$
Otherwise, there is $R>0$ such that
$$
\lim_{n \to +\infty}\sup_{z \in \mathbb{R}^N}\int_{B_R(z)}|u_n|^{2}\,dx=0.
$$ 
Hence, by Lemma \ref{lions}, 
$$
\lim_{n \to +\infty}\int_{\mathbb{R}^N}|u_n|^q\,dx=0, \quad \ q\in (2, 2^{*}_{s}).
$$
Since  $0<c< \frac{s}{N} (2^{-1}C(N,s)S)^{N/2s}$, arguing as in \cite[Lemma 3.3]{carrion} for local case, \cite[Lemma 3.4]{Shang} for nonlocal case, the last limit combined with $I'_\epsilon(u_n)u_n=o_n(1)$ gives
$$
u_n \to 0 \quad \mbox{in} \quad H^{s}(\mathbb{R}^N),
$$ 
implying that
$$
I_\epsilon(u_n) \to 0,
$$
which is a contradiction, because by hypotheses $I_\epsilon(u_n) \to c>0$. Thereby, for each $R>0$,  
there are $(z_n) \subset \mathbb{R}^N$, $\tau>0$ and a subsequence of $(u_n)$, still denoted by itself, such that
\begin{equation} \label{L 1}
\int_{B_R(z_n)}|u_n|^{2}\,dx \geq \tau, \quad \forall n \in \mathbb{N}.
\end{equation} 
Setting $v_n=u_n(\cdot+z_n)$, we have that $(v_n)$ is bounded in $H^{s}(\mathbb{R}^N)$. Thus, for 
some subsequence of $(u_n)$, still denoted by $(u_n)$, 
there is $v \in H^{s}(\mathbb{R}^N)$ such that
\begin{equation} \label{L2}
v_n \rightharpoonup v \quad \mbox{in} \quad H^{s}(\mathbb{R}^N).
\end{equation} \label{L3}
From (\ref{L 1}) and (\ref{L2}),
$$
\int_{B_R(0)}|v|^{2}\,dx \geq \tau,
$$
showing that $v \not= 0$. Moreover, (\ref{L 1}) gives that $|z_n| \to +\infty$, because 
$u_n \rightharpoonup 0$ in $H^{s}(\mathbb{R}^N)$.
\fim

\vspace{0.5 cm}

The lemma below brings an important estimate from above involving the mountain pass level $m_\lambda(c_0)$, which is crucial in our approach.

\begin{lem} \label{ESTIMATIVA SUPERIOR} There is $\lambda^{*}>0$ such that 
$$
m_\lambda(c_0) \leq \frac{s}{2N} (2^{-1}C(N,s)S)^{N/2s},
$$
for all $\lambda \geq \lambda^{*}$.	
\end{lem}
\noindent {\bf Proof.}
Let  $w \in H^{s}(\mathbb{R}^{2}) \setminus \{0\}$. We know that there is $t_\lambda>0$ such that 
$$
t_\lambda w \in \mathcal{M}_{c_0}=\left\{u \in H^{s}(\mathbb{R}^{N}) \setminus \{0\}\,:\, J'_{\lambda,c_0}(u)u=0\right\},
$$ 
that is
$$
\|w\|^{2}=\lambda t_\lambda^{q-2} |w|_{q}^{q}+t_\lambda^{2_s^{*}-2}|w|_{2_s^{*}}^{2_s^{*}}.
$$ 
The above equality gives
$$
t_\lambda \to 0 \quad \mbox{as} \quad \lambda \to +\infty.
$$
As
$$
m_\lambda(c_0) \leq \max_{t \geq 0}J_{\lambda,c_0}(tw)\leq \frac{t_{\lambda}^{2}}{2}\|w\|^{2},
$$
we derive that
$$
m_\lambda(c_0) \to 0 \quad \mbox{as} \quad \lambda \to +\infty, 
$$
finishing the proof. \fim

\vspace{0.5 cm}

As a consequence of the last results we have the following corollary

{\begin{cor} \label{ground state} The problem 
$$
(-\Delta)^{s}{u}+c_0u=\lambda |u|^{q-2}u+|u|^{2^{*}_s-2}u\,\,\, \mbox{in} \,\,\, \mathbb{R}^{N}, \eqno{(P_{\lambda,\infty})}
$$
possesses a positive ground state solution for all $\lambda \geq \lambda^{*}$, that is, there is $w \in H^{s}(\mathbb{R}^{N})$ such that
$$
J_{\lambda,c_0}(w)=m_\lambda(c_0) \quad \mbox{and} \quad J'_{\lambda,c_0}(w)=0,
$$	
for all $\lambda \geq \lambda^{*}$.
\end{cor}
\noindent {\bf Proof.}  The existence of a ground state solution can be obtained repeating the same idea found in Alves, Carri\~ao and Miyagaki \cite{carrion} for the local case, that is, $s=1$. By using the definition $\|\,\,\,\,\|$, we have that
$$
\||u|\| \leq \|u\| \quad \forall u \in H^{s}(\mathbb{R}^{N}).
$$ 
Thus, if $u_0$ is a ground state solution,   
$$
J'_{\lambda,c_0}(|u_0|)|u_0| \leq J'_{\lambda,c_0}(u_0)u_0. 
$$
From this, there is $t_1 \in (0,1]$ such that $t_1 |u_{0}|\in \mathcal{M}_{c_0}$, and so, 
$$
m_\lambda(c_0) \leq J_{\lambda,c_0}(t_1|u_{0}|) \leq J_{\lambda,c_0}(u_0)=m_\lambda(c_0), 
$$
implying that $J_{\lambda,c_0}(t_1|u_{0}|)=m(c_0)$. Using Deformation Lemma, we deduce that $t_1 |u_{0}|$ is a critical point, and so, it is a ground state solution, finishing the proof.
\fim

\vspace{0.5 cm}
The next lemma shows that there is positive radial ground state solution.
\begin{lem} \label{RADIAL} If $u_0$ is a positive ground state solution of  $(P_{\lambda,\infty})$, then its symmetrization denoted by $u^{*}_0$ is also a positive ground state solution of  $(P_{\lambda,\infty})$.  
\end{lem} 
\noindent {\bf Proof.} Denote by $u^{*}_0$ the symmetrization of $u_0$. Then, c.f. \cite{Park},
$$
\|u^{*}_{0}\|\leq \|u_0\|, \quad |u^{*}_{0}|_q=|u_{0}|_q \quad \mbox{and} \quad |u^{*}_{0}|_{2^{*}_{s}}=|u_{0}|_{2^{*}_{s}}.
$$ 

From this, $J'_{\lambda,c_0}(u^{*}_{0})u^{*}_{0}\leq 0$. Then, there is $t_0 \in (0,1]$ such that $t_0 u^{*}_{0} \in \mathcal{M}_{c_0}$. Thereby,
$$
m_\lambda(c_0) \leq J_{\lambda,c_0}(t_0u^{*}_{0}) \leq J_{\lambda,c_0}(u^{*})=m_\lambda(c_0), 
$$
implying that $J_{\lambda,c_0}(t_0u^{*}_{0})=m_\lambda(c_0)$. The above inequality also ensures that $t_0=1$, otherwise we must have
$$
m_\lambda(c_0) \leq J_{\lambda,c_0}(t_0u^{*}_{0}) < J_{\lambda,c_0}(u^{*})=m_\lambda(c_0), 
$$
which is an absurd. From this, $J_{\lambda,c_0}(u^{*}_{0})=m_\lambda(c_0)$. Now, applying  Deformation Lemma, we deduce that $u^{*}_{0}$ is critical point, then it is a radial ground state solution. \fim

\vspace{0.5 cm}

The lemma below is a key point in our arguments, because it is a regularity result for problems of the type
$$
\left\{
\begin{array}{l}
(-\Delta)^{s}{u}+\alpha u=\lambda |u|^{q-2}u+|u|^{2^{*}_s-2}u,\,\,\, \mbox{in} \,\,\, \mathbb{R}^{N}, \\
u \in H^{s}(\mathbb{R}^{N})
\end{array}
\right.
\eqno{(P_{\lambda,\alpha})}
$$
for $\lambda, \alpha>0$
\begin{lem} \label{regularidade} If $u \in H^{s}(\mathbb{R}^{N})$ is a solution of $(P_{\lambda,\alpha})$, then $u \in C^{2}(\mathbb{R}^{N}) \cap H^{1}(\mathbb{R}^{N})$.	
\end{lem}
\noindent {\bf Proof.} In what follows, we will use an approach due to Cabr\'e  and Sire \cite{cabre} , that is, we will see the problem of the following way , c.f. \cite{caffarelli},
$$
\left\{
\begin{array}{rclcl}
div(y^{1-2s}{\nabla v})&=&0,  & \mbox{in} & \mathbb{R}_{+}^{N+1}, \\
2(1-s)\frac{\partial v}{\partial \nu^{s}}&=&-\alpha v+f(v), & \mbox{on}& \mathbb{R}^{N},
\end{array}
\right.
\eqno{(P^{*}_{\lambda,\alpha})}
$$
where $\mathbb{R}_{+}^{N+1}=\{(x_1,....,x_N,y) \in \mathbb{R}^{N+1}\,:\, y>0\}$, $\lambda, \alpha>0$ and 
$$
\frac{\partial v}{\partial \nu^{s}}(x)=- \lim_{y \to 0^{+}}y^{1-2s}\frac{\partial v}{\partial y }(x,y).
$$

Associated with $(P^{*}_{\lambda,\alpha})$, we have the energy functional $I:X^{1,\alpha} \to \mathbb{R}$ given by
\begin{equation} \label{E1}
I(v)=\frac{1}{2}\int_{\mathbb{R}_{+}^{N+1}}y^{1-2s}|\nabla v|^{2}\,dxdy+\frac{1}{2}\int_{\mathbb{R}^{N}}\alpha|v|^{2}\,dx-\int_{\mathbb{R}^{N}}F(v)\,dx
\end{equation}
where  $F$ denotes the primitive of $f$, that is,
$$
F(t)=\int_{0}^{t}f(s)\,ds=\lambda \frac{1}{q}|t|^q +\frac{1}{2^{*}_{s}}|t|^{2^{*}_{s}}, \quad \forall t \in \Re
$$
and $X^{1,s}$ is the Hilbert space obtained as the closure of $C_{0}^{\infty}(\overline{\mathbb{R}_{+}^{N+1}})$ under the norm
$$
\|v\|_{1,s}=\left( \int_{\mathbb{R}_{+}^{N+1}}|\nabla v|^{2}\,dxdy + \int_{\mathbb{R}^{N}}\alpha|v|^{2}\,dx\right)^{\frac{1}{2}}.
$$
Using some embeddings mentioned in  Br\"andle,Colorado and  S\'anchez \cite{Brandle}( see also \cite{CW,ZhangLiuJiao}), we deduce that the embeddings
$$
X^{1,s} \hookrightarrow L^{p}(\mathbb{R}^{N}) \,\,\, \mbox{for} \,\,\ p \in [2, 2^{*}_{\alpha}] 
$$
are continuous, where $2^{*}_{s}=\frac{2N}{N-2s}$. Moreover, we know that $u$ is a solution of $(P_{\lambda,\alpha})$ if, and only if, $u=v(x,0)$ for all $x \in \mathbb{R}^{N}$, for some critical point $v$ of $I$.

In what follows,  for each $L>0$, we set  
$$
v_L(x,y)=
\left\{
\begin{array}{lcr} 
v(x,y), &\mbox{if}& (x,y)\leq L\\
L, &\mbox{if}& v(x,y) \geq L
\end{array}
\right.
$$
and
$$
z_L=v_L^{2(\beta -1)}v,
$$
with $\beta >1$ to be determined later. Since $I'(v)z_L=0,$  adapting the same approach explored in Alves and Figueiredo \cite[Lemma 4.1]{AF2}, we will find the following estimate
$$
|v(.,0)|_{\infty} \leq C |v(.,0)|_{{2^{*}_s}},
$$
or equivalently, 
$$
|u|_{{\infty}} \leq C |u|_{{2^{*}_s}}.
$$

Now, fixing $M=|u|_\infty+1$, we consider the following function
$$
g_M(t)=
\left\{
\begin{array}{lcc}
0, & \mbox{if} & t \leq 0 \\
\lambda t^{q-1}+t^{2^{*}_s}-1, & \mbox{if} & 0 \geq t \leq M, \\
\lambda t^{q-1}+ A_Mt^{q-1}+B_M, & \mbox{if} & t \geq M,
\end{array}
\right.
$$
where $A_M$ and $B_M$ are chosen  such that $g_M \in C^{1}(\mathbb{R})$. It is easy to see that $g_M$ has a subcritical growth and $u$ is a solution of the problem
$$
\left\{
\begin{array}{l}
(-\Delta)^{s}{u}+\alpha u=g_M(u),\,\,\, \mbox{in} \,\,\, \mathbb{R}^{N}, \\
u \in H^{s}(\mathbb{R}^{N}).
\end{array}
\right.
\eqno{(P_{\lambda,\alpha,M})}
$$
Using the arguments explored in Felmer, Quass and Tan \cite{FQT}( see also \cite{Frank1}), we deduce that
$$
|u(x)| \to 0 \quad \mbox{as} \quad |x| \to +\infty.  
$$

This way, we see that $(-\Delta)^{s}u \in L^{\infty}(\mathbb{R}^{N})$,  and so, $u \in C^{2}(\mathbb{R}^{N})$. Repeating the same arguments found in \cite[Section 2]{Moustapha}, we also have $|\nabla u| \in L^{2}(\mathbb{R}^{N})$. As $u \in L^{2}(\mathbb{R}^{N})$, it follows that $u \in H^{1}(\mathbb{R}^{N})$.

\fim

\begin{lem} \label{L1} Under the hypotheses $(V_1)-(V_4)$ and $\lambda \geq \lambda^{*}$, for each $\sigma >0$,
 there is $\epsilon_0=\epsilon(\lambda, \sigma)>0$, such that $I_\epsilon$ satisfies the $(PS)_c$ condition for
 all $c \in (m_\lambda(c_0)+\sigma,2m_\lambda(c_0)-\sigma)$, for all $\epsilon \in (0, \epsilon_0)$.  
\end{lem}

\noindent {\bf Proof.} \, We will prove the lemma arguing by contradiction. Suppose that there is $\sigma >0$ and 
$\epsilon_n \to 0$, such that $I_{\epsilon_n}$ does not satisfy the $(PS)$ condition.  

Thereby, there is $c_n \in  (m_\lambda(c_0)+\sigma,2m_\lambda(c_0)-\sigma)$ such that
  $I_{\epsilon_n}$ does not verify the $(PS)_{c_n}$ condition. Then, there  is sequence $(u^{n}_m)$ such that
\begin{equation} \label{E1}
\lim_{m \to +\infty}I_{\epsilon_n}(u^{n}_m)=c_n \quad \mbox{and} \quad \lim_{m \to +\infty}I'_{\epsilon_n}(u^{n}_m)=0,
\end{equation}
with
\begin{equation} \label{E2}
 u^{n}_m \rightharpoonup u_n \quad \mbox{in} \quad H^{s}(\mathbb{R}^N) \quad \mbox{but} 
\quad u^{n}_m \not\to u_n \quad \mbox{in} \quad H^{s}(\mathbb{R}^N).
\end{equation}

Then, for $v^{n}_m=u^{n}_m-u_n$, the Brezis-Lieb Lemma  yields 
$$
I_{\epsilon_n}(u^{n}_m)=I_{\epsilon_n}(u_n)+I_{\epsilon_n}(v^{n}_m)+o_m(1) \quad
 \mbox{and} \quad I'_{\epsilon_n}(v^{n}_m)=o_m(1).
$$

\begin{claim} \label{CorLions} There is $\delta>0$, such that
$$
\liminf_{m \to +\infty}\sup_{y \in \mathbb{R}^N}\int_{B_R(y)}|v^{n}_{m}|^{2}\,dx \geq \delta, \quad \forall n \in \mathbb{N}.
$$	
\end{claim}
Indeed, if the claim does not hold, there is $(n_j) \subset \mathbb{N}$ satisfying
$$
\liminf_{m \to +\infty}\sup_{y \in \mathbb{R}^N}\int_{B_R(y)}|v^{n_j}_{m}|^{2}\,dx \leq \frac{1}{j}, \quad \forall j \in \mathbb{N}.
$$
Using the arguments found in \cite{W} and \cite{Secchi}, we deduce that
\begin{equation} \label{Lq}
\limsup_{m \to +\infty}|v^{n_j}_{m}|_q=o_j(1), \quad \forall q \in (2, 2^{*}_{s}).
\end{equation}
Then
$$
\limsup_{m \to +\infty}\int_{\mathbb{R}^N}|v^{n_j}_{m}|^q\,dx=o_j(1).
$$
As $I_{\epsilon_n}(u_n)\geq 0$ and $c <\frac{s}{N} (2^{-1}C(N,s)S)^{N/2s},$ the above estimate combined with (\ref{Lq})  and  $I'_{\epsilon_{n_j}}(v^{n_j}_{m})(v^{n_j}_{m})=o_{m}(1)$ gives
$$
\limsup_{m \to +\infty}\|v^{n_j}_m\|^{2}=o_j(1).
$$
Now since  $u^{n_j}_m \not\to u_{n_j}$ in $H^{s}(\mathbb{R}^N)$, we derive that
$$
\liminf_{m \to +\infty}\|v^{n_j}_m\|^{2}>0.
$$
Then, without loss of generality, we can assume that $(v^{n_j}_m) \subset H^{s}(\mathbb{R}^N) \setminus \{0\}$. 
Thereby, there is $t^{n_j}_m \in (0,+\infty)$ such that
$$
t^{n_j}_m v^{n_j}_m \in \mathcal{N}_{\epsilon_{n_j}},
$$
\noindent verifying
$$
\lim_{m \to +\infty}t^{n_{j}}_m=1
\quad \mbox{
and}\quad
\lim_{m \to +\infty}I_{\epsilon_{n_j}}(t^{n_{j}}_m v^{n_j}_m)= \lim_{m \to +\infty}I_{\epsilon_{n_j}}(v^{n_j}_m).
$$ 
From the above informations, there is $r^{n_j}_m \in (0,1)$ such that
$$
r^{n_j}_m (t^{n_j}_m v^{n_j}_m) \in \mathcal{M}_{c_0}.
$$
Hence,
\begin{eqnarray*}
m_\lambda(c_0) &\leq& \limsup_{m \to +\infty}J_{c_0}(r^{n_j}_m (t^{n_j}_m v^{n_j}_m))
\leq \limsup_{m \to +\infty}I_{\epsilon_{n_j}}(t^{n_j}_m v^{n_j}_m)\\ & =&
\limsup_{m \to +\infty}I_{\epsilon_{n_j}}(v^{n_j}_m)\leq \frac{(1+|V|_{\infty})}{2} \limsup_{m \to +\infty}\|v^{n_j}_m\|^{2},
\end{eqnarray*}
that is,
$$
m_\lambda(c_0) \leq o_j(1),
$$
which is a contradiction.

From the above study, for each $m \in \mathbb{N}$, there is $m_n \in \mathbb{N}$ such that
$$
\int_{B_R(z^{n}_{m_n})}|u^{n}_{m_n}|^{2}\,dx \geq \frac{\delta}{2}, \quad |\epsilon_n z^{n}_{m_n}| \geq n, 
\quad \|I'_{\epsilon_n}(u^{n}_{m_n})\|\leq \frac{1}{n} \quad \mbox{and} \quad |I_{\epsilon_n}(u^{n}_{m_n})-c_n|\leq \frac{1}{n}.
$$

In what follows, we denote by $(z_n)$ and $(u_n)$ the sequences $(z^{n}_{m_n})$ and $(u^ {n}_{m_n})$ respectively.  Then,
$$
\int_{B_R(z_n)}|u_{n}|^{2}\,dx \geq \frac{\delta}{2}, \quad |\epsilon_n z_{n}| \geq n, \quad 
\|I'_{\epsilon_n}(u_n)\|\leq \frac{1}{n} \quad \mbox{and} \quad |I_{\epsilon_n}(u_n)-c_n|\leq \frac{1}{n}.
$$

\begin{claim} \label{limite fraco}  $u_n \rightharpoonup 0$ in $H^{s}(\mathbb{R}^N)$. 
	
\end{claim}

Indeed, assume by contradiction that there is $u \in H^{s}(\mathbb{R}^N) \setminus \{0\}$ such that
$$
u_n \rightharpoonup u \quad \mbox{in} \quad H^{s}(\mathbb{R}^N). 
$$
Using the limit $ \|I'_{\epsilon_n}(u_n)\|\to 0$, it is possible to prove that $u$ is a solution of the problem
$$
(-\Delta)^{s}{u}+V(0)u-f(u)=0 \quad \mbox{in} \quad \mathbb{R}^N. 
$$

Then,  the definition of $m_\lambda(V(0))$ together with  $(V_4)$ gives  
$$
J_{\lambda, V(0)}(u) \geq m_\lambda(V(0)) \geq 2m_\lambda(c_0).
$$
On the other hand,  the Fatous' lemma leads to
$$
J_{\lambda, V(0)}(u) \leq \liminf_{n \to +\infty}[I_{\epsilon_n}(u_n)-\frac{1}{\theta}I'_{\epsilon_n}(u_n)]=
\liminf_{n \to +\infty}I_{\epsilon_n}(u_n)=\liminf_{n \to +\infty}c_n \leq 2m_\lambda(c_0)-\sigma,
$$
obtaining a contradiction. Then , the Claim \ref{limite fraco} is proved.

Considering $w_n=u_n(\cdot+z_n)$, we have that $(w_n)$ is bounded in $H^{s}(\mathbb{R}^N)$. 
Then, there is $w \in H^{s}(\mathbb{R}^N)$ such that
$$
w_n \rightharpoonup w \quad \mbox{in} \quad H^{s}(\mathbb{R}^N). 
$$
Hence,
$$
\int_{B_R(0)}|w|^{2}\,dx \geq \frac{\delta}{2},
$$
showing that $w \not= 0$.

Now, for each $\phi \in H^{s}(\mathbb{R}^N)$, we have the equality below
$$
\int_{\mathbb{R}^N}|\xi|^{2s} \widehat{ w_n}\widehat{ \phi} \, d \xi + \int_{\mathbb{R}^N}V(\epsilon_n z_n+\epsilon_n z)w_n \phi
 \, dx - \int_{\mathbb{R}^N}f(w_n)\phi \, dx = o_n(1)\|\phi\|
$$
which implies that, see \cite[Theorem 3.5]{Shang},  $w$ is a nontrivial solution of the problem
\begin{equation} \label{equacao}
(-\Delta)^{s}{u}+\alpha_1u-f(u)=0 \quad \mbox{in} \quad \mathbb{R}^N,
\end{equation}
where $\alpha_1=\displaystyle \lim_{n \to +\infty}V(\epsilon_n z_n)$. Thereby, by Lemma \ref{regularidade}, $w \in C^{2}(\mathbb{R}^N) \cap H^{1}(\mathbb{R}^N)$.

For each $k \in \mathbb{N}$, there is $\phi_k \in C^{\infty}_{0}(\mathbb{R}^{N})$ such that
$$
\|\phi_k -w\| \to 0 \quad \mbox{as} \quad k \to +\infty,
$$
that is,
$$
\|\phi_k -w\|=o_k(1).
$$

Using $\frac{\partial \phi_k}{\partial x_i}$ as a test function, we get
$$
\int_{\mathbb{R}^N}|\xi|^{2s} \widehat{ w_n}\widehat{  \frac{\partial \phi_k}{\partial x_i}} d \xi+  
\int_{\mathbb{R}^N}V(\epsilon_n z+\epsilon_n z_n)w_n\frac{\partial \phi_k}{\partial x_i}\,dx-
\int_{\mathbb{R}^N}f(w_n)\frac{\partial \phi_k}{\partial x_i}\,dx=o_n(1).
$$
Now, using well known arguments, we have that
$$
\int_{\mathbb{R}^N}|\xi|^{2s}  \widehat{ w_n}\widehat{  \frac{\partial \phi_k}{\partial x_i}} d \xi=\int_{\mathbb{R}^N}|\xi|^{2s}  \widehat{ w}\widehat{  \frac{\partial \phi_k}{\partial x_i}} d \xi+o_n(1)
$$
and
$$
\int_{\mathbb{R}^N}f(w_n)\frac{\partial \phi_k}{\partial x_i}\,dx=\int_{\mathbb{R}^N}f(w)\frac{\partial \phi_k}{\partial x_i}\,dx+o_n(1).
$$
Gathering the above limit with (\ref{equacao}), we deduce that  
$$
\limsup_{n \to +\infty}\left|\int_{\mathbb{R}^N}(V(\epsilon_n z_n+\epsilon_n z)-V(\epsilon_n z_n))w_n\frac{\partial \phi_k }
{\partial x_i}\,dx\right|=0.
$$
As $\phi_k$ has compact support , the above limit gives 
$$
\limsup_{n \to +\infty}\left|\int_{\mathbb{R}^N}(V(\epsilon_n z_n+\epsilon_n z)-V(\epsilon_n z_n))w\frac{\partial \phi_k}{\partial x_i}\,dx\right|=0.
$$
Also since  $\frac{\partial w}{\partial x_i} \in L^{2}(\mathbb{R}^N)$, we have that $(\frac{\partial \phi_k}{\partial x_i})$ 
is bounded in $L^{2}(\mathbb{R}^N).$ Hence, 
$$
\limsup_{n \to +\infty}\left|\int_{\mathbb{R}^N}(V(\epsilon_n z_n+\epsilon_n z)-V(\epsilon_n z_n))
\phi_k\frac{\partial \phi_k}{\partial x_i}\,dx\right|=o_k(1),
$$
and so,
$$
\limsup_{n \to +\infty}\left|\frac{1}{2}\int_{\mathbb{R}^N}(V(\epsilon_n z_n+\epsilon_n z)-V(\epsilon_n z_n))\frac{\partial 
(\phi_k^{2})}{\partial x_i}\,dx\right|=o_k(1).
$$
Using Green's Theorem together with the fact that $\phi_k$ has compact support, we find the limit below
$$
\limsup_{n \to +\infty}\left|\int_{\mathbb{R}^N}\frac{\partial V}{\partial x_i}(\epsilon_n z_n+\epsilon_n z) 
\, \phi_k^{2}\,dx\right|=o_k(1),
$$
which leads to 
$$
\limsup_{n \to +\infty}\left|\frac{\partial V}{\partial x_i}(\epsilon_n z_n)\int_{\mathbb{R}^N}|\phi_k|^{2}\,dx\right|=o_k(1).
$$
As 
$$
\int_{\mathbb{R}^N}|\phi_k|^{2}\,dx \to \int_{\mathbb{R}^N}|w|^{2}\,dx \quad \mbox{as} \quad k \to +\infty,
$$
it follows that
$$
\limsup_{n \to +\infty}\left|\frac{\partial V}{\partial x_i}(\epsilon_n z_n)\right|=o_k(1), \quad \forall i \in \{1,....,N\}.
$$
Since $k$ is arbitrary, we derive that 
$$
\nabla V(\epsilon_n z_n) \to 0 \quad \mbox{as} \quad n \to \infty.
$$
Therefore,  $(\epsilon_n z_n)$ is a $(PS)_\alpha$ sequence for $V$, which is an absurd, because by hypotheses on  $V$, it  satisfies the $(PS)$ condition and $(\epsilon_n z_n)$ does not have any convergent subsequence in $\mathbb{R}^N$. 
\fim

\vspace{0.5 cm}

Denote by ${\mathcal N}_{\epsilon}$ the Nehari Manifold associated with $I_\epsilon$, that is, 
$$ 
 {\mathcal N}_{\epsilon}=\left\{ u \in H^{s}(\mathbb{R}^N) \setminus \{0\}\,:\, I'_{\epsilon}(u)u=0    \right\}.
$$

\begin{lem}\label{PS em Nehari} For $\lambda \geq \lambda^{*}$ and $\sigma>0$, the  functional $I_\epsilon$ restrict to ${\mathcal N}_{\epsilon}$  satisfies the $(PS)_c$ condition for all $c \in (m_\lambda(c_0)+\sigma,2m_\lambda(c_0)-\sigma)$. 
	\end{lem}
\noindent {\bf Proof.} Let $(u_n)$ be a $(PS)$-sequence for
$I_{\epsilon}$ constrained to $\mathcal{M}_{\epsilon}$. Then
$I_{\epsilon}(u_{n})\rightarrow c$ and
\begin{eqnarray}\label{contra}
I'_{\epsilon}(u_{n}) = \theta_{n} G_{\epsilon}'(u_{n}) + o_{n}(1),
\end{eqnarray}
for some $(\theta_{n}) \subset \mathbb{R}$, where
$G_{\epsilon}:H^{s}(\mathbb{R}^N) \rightarrow \mathbb{R}^N$ is
given by
\begin{eqnarray*}
	G_{\epsilon}(v) := \displaystyle\int_{\mathbb{R}^N}|\xi|^{2s} |\widehat{ v}|^{2}d \xi + 
\int_{\mathbb{R}^N}V(\epsilon x)|v|^{2} \,dx
	-\displaystyle\int_{\mathbb{R}^N}f(v)v \,dx.
\end{eqnarray*}
Notice that $G'_{\epsilon}(u_{n})u_{n}\leq 0$. By standard
arguments show that $(u_n)$ is bounded. Thus, up to a subsequence,
$G_{\epsilon}'(u_{n}) u_{n} \rightarrow l\leq 0$. If $l \neq 0$, we infer
from (\ref{contra}) that $\theta_{n}=o_{n}(1)$. In this case, we can
use (\ref{contra}) again to conclude that $(u_{n})$ is a $(PS)_{c}$
sequence for $I_{\epsilon}$ in $H^{s}(\mathbb{R}^N)$, and so, $(u_{n})$ has a strongly convergent subsequence. If
$l=0$, it follows that 
$$ 
\displaystyle\int_{\mathbb{R}^N}(f'(u_n)u_n^{2}-f(u_n)u_n) \ dx\rightarrow 0.
$$
Using the definition of $f$, we know that
\begin{equation} \label{Desig}
f'(t)t^{2}-f(t)t>0, \quad \forall t \in \mathbb{R} \setminus \{0\}.
\end{equation}
If $u \in H^{s}(\mathbb{R}^N)$ is the weak limit of $(u_n)$, the Fatous' Lemma combined with the last limit leads to
$$
\int_{\mathbb{R}^N}(f'(u)u^{2}-f(u)u) \ dx=0.
$$
Then, by (\ref{Desig}), $u=0$. Applying Corollary \ref{sequencia}, there is $(y_n)
 \subset \mathbb{R}^N$ with $|y_n| \to +\infty$ such that
$$
v_n=u_n(\cdot +y_n) \rightharpoonup v \not= 0 \quad \mbox{in} \quad H^{s}(\mathbb{R}^N).
$$
By change variable, 
$$
\displaystyle\int_{\mathbb{R}^N}(f'(v_n)v_n^{2}-f(v_n)v_n)=\displaystyle\int_{\mathbb{R}^N}(f'(u_n)u_n^{2}-f(u_n)u_n)
 \ dx\rightarrow 0.
$$
Applying again Fatous's Lemma, we get
$$
\int_{\mathbb{R}^N}(f'(v)v^{2}-f(v)v) \ dx=0,
$$
which is an absurd, because being  $v \not=0$, the inequality (\ref{Desig}) leads to
$$
\int_{\mathbb{R}^N}(f'(v)v^{2}-f(v)v) \ dx>0,
$$
finishing the proof of the lemma.

\fim

\begin{cor} \label{ponto critico}
If $u \in H^{s}(\mathbb{R}^N)$ is a critical point of $I_\epsilon$ restrict to ${\mathcal N}_{\epsilon}$,
 then $u$ is a critical point of  $I_\epsilon$ in $H^{s}(\mathbb{R}^N)$. 
\end{cor}
\noindent {\bf Proof.} The proof follows arguing as in the proof of Lemma \ref{PS em Nehari}. \fim

\vspace{0.5 cm}

The next  lemma will be crucial in our study to show a lower estimate involving  a special minimax level,  
which will be defined later on. 
\begin{lem} \label{compacidade} Let $\epsilon_n \to 0$ and $(u_n) \subset \mathcal{N}_{\epsilon_n}$ such
 that $I_{\epsilon_n}(u_n) \to m_\lambda(c_0)$. Then, there is $(z_n) \subset \mathbb{R}^N$ with
 $|z_n| \to +\infty$ and $u_1 \in H^{s}(\mathbb{R}^N) \setminus \{0\}$ such that
$$
u_n(\cdot +z_n) \to u_1 \quad \mbox{in} \quad H^{s}(\mathbb{R}^N).
$$
Moreover, $\displaystyle\liminf_{n \to +\infty}|\epsilon_n z_n|>0$.	
\end{lem}
\noindent {\bf Proof.} Since $u_n \in \mathcal{N}_{\epsilon_n}$, we have that
$J_{\lambda,c_0}'(u_n)u_n<0$ for all $n \in \mathbb{N}$. Thus, there is $t_n \in (0, 1)$ such that 
$t_n u_n \in \mathcal{M}_{c_0}.$ Therefore, 
$$
(t_n u_n) \subset  \mathcal{M}_{c_0} \quad \mbox{and} \quad J_{\lambda,c_0}(t_n u_n)\to m_\lambda(c_0).
$$
Now, by \cite[Lemma 5.1]{Shang}, there are $(z_n) \subset \mathbb{R}^N$, $u_1 \in H^{s}(\mathbb{R}^N) 
\setminus \{0\}$, and a subsequence of $(u_n)$, still denoted by $(u_n)$,  verifying
$$
u_n(\cdot+z_n) \to u_1 \quad \mbox{in} \quad H^{s}(\mathbb{R}^N).
$$
\begin{claim} \label{zn} 
$\displaystyle \liminf_{n \to +\infty}|\epsilon_n z_n|>0 $.	
\end{claim}
Indeed, as $u_n \in \mathcal{N}_{\epsilon_n}$ for all $n \in \mathbb{N}$, the function $u_n^{1}=u_n(\cdot+z_n)$ must verify
\begin{equation} \label{Eq1}
\int_{\mathbb{R}^N}|\xi|^{2s}|\widehat{ u_{n}^{1}}|^{2}  d \xi+ \int_{\mathbb{R}^N} V(\epsilon_n x+\epsilon_n z_n)|u_n^{1}|^{2}
\,dx=\int_{\mathbb{R}^N}f(u_n^{1})u_n^{1}\,dx.
\end{equation}
Supposing by contradiction, up to a  subsequence, 
$$
\lim_{n \to +\infty}\epsilon_n z_n=0.
$$
Taking the limit of $n \to +\infty$ in (\ref{Eq1}), we get
$$
\int_{\mathbb{R}^N}|\xi|^{2s}|\widehat{ u_{1}}|^{2}  d \xi+ \int_{\mathbb{R}^N} V(0) |u_1|^2 dx =
\int_{\mathbb{R}^N}f(u_{1})u_{1}\,dx.
$$
Thereby, 
$$
J_{\lambda, V(0)}(u_{1}) \geq m_\lambda(V(0)) > m_\lambda(c_0).
$$
On the other hand, 
$$
I_{\epsilon_n}(u_n) \to J_{\lambda, V(0)}(u_{1}),
$$ 
which leads to
$$
m_\lambda(c_0)=J_{\lambda, V(0)}(u_{1}),
$$
obtaining a contradiction. 
\fim	
	
\section{A special minimax level}

In order to prove the Theorem \ref{T2}, we will consider  a special minimax level. The construction involves the barycenter function used in \cite{Benci}, see also \cite{AlvesNovo,PFM}, given by
$$
\beta(u)= \frac{\displaystyle  \int_{\mathbb{R}^N}\frac{x}{|x|}|u|^{2}\,dx}{\displaystyle \int_{\mathbb{R}^N}|u|^{2}\,dx},
 \quad \forall u \in H^{s}(\mathbb{R}^N) \setminus \{0\}.
$$

\noindent For each $z \in \mathbb{R}$ and $\epsilon >0$, let us define the function
$$
\phi_{\epsilon, z}(x)=t_{\epsilon,z}u_0\left( x - \frac{z}{\epsilon} \right),
$$
where $t_{\epsilon,z}>0$ is such that $\phi_{\epsilon, z} \in {\mathcal N}_{\epsilon}$, and  $u_0 \in H^{s}(\mathbb{R}^{N})$ 
is a radial positive ground state solution $u_0 \in H^{s}(\mathbb{R}^{N})$ for $J_{\lambda, c_0}$, that is, 
$$
J_{\lambda, c_0}(u_0)=m_\lambda(c_0) \quad \mbox{and} \quad J'_{\lambda, c_0}(u_0)=0,
$$ 
\noindent whose the existence was guaranteed in Lemma \ref{RADIAL}.

\vspace{0.5 cm}

We establish several properties involving $\beta$ and $\phi_{\epsilon, z}.$

\begin{lem} \label{B1} 
	For each  $r>0$, $\displaystyle \lim_{\epsilon \to 0}\left(\sup\left\{\left|\beta(\phi_{\epsilon,z})-\frac{z}{|z|} 
\right|\,:\,|z| \geq r \right\}\right)=0$.
\end{lem}
\noindent {\bf Proof.} It is enough to show  that for any $(z_n)$ with $|z_n| \geq r$ and $\epsilon_n \to 0$, we have that 
$$
\left|\beta(\phi_{\epsilon_n,z_n})-\frac{z_n}{|z_n|} \right| \to 0 \quad \mbox{as} \quad n \to +\infty.
$$
By change variable, 
$$
\left|\beta(\phi_{\epsilon_n,z_n})-\frac{z_n}{|z_n|} \right|=\frac{\displaystyle \int_{\mathbb{R}^N}\left|\frac{\epsilon_n x 
+z_n}{|\epsilon_n x +z_n|}-\frac{z_n}{|z_n|}\right||u_0(x)|^{2}\,dx}{\displaystyle \int_{\mathbb{R}^N}|u_0|^{2}\, dx}.
$$
Since for each $x \in \mathbb{R}$,
$$
\left|\frac{\epsilon_n x +z_n}{|\epsilon_n x +z_n|}-\frac{z_n}{|z_n|}\right| \to 0 \quad \mbox{as} \quad n \to +\infty, 
$$
applying the  Lebesgue dominated convergence Theorem, we get
$$
\int_{\mathbb{R}^N}\left| \frac{\epsilon_n x +z_n}{|\epsilon_n x +z_n|}-\frac{z_n}{|z_n|}\right||u_0(x)|^{2}\,dx \to 0. 
$$
This proves  the lemma.  \fim

\vspace{0.5 cm}

As an immediate consequence, we have 
\begin{cor} \label{cor1} 
	Fixed $r>0$, there is $\epsilon_0>0$ such that
	$$
	(\beta(\phi_{\epsilon,z}),z)>0, \quad \forall |z| \geq r \quad \mbox{and} \quad \forall \epsilon \in (0,\epsilon_0).
	$$
\end{cor}

\vspace{0.5 cm}

Define now the set 
$$
{\mathcal B}_{\epsilon}=\{u \in {\mathcal N}_{\epsilon}\,:\, \beta(u) \in Y  \}.
$$
Note that ${\mathcal B}_{\epsilon} \not=\emptyset$, because $\phi_{\epsilon,0}=0 \in Y,$ for all $\epsilon>0$. 
Associated with  the above set, define  the real number $D{_\epsilon}$ given by
$$
D_{\epsilon}=\inf_{u \in \mathcal{B}_{\epsilon}}I_{\epsilon}(u).
$$

\vspace{0.5 cm}

Next lemma establishes an important relation between  $D_{\epsilon}$ and $m_\lambda(c_0)$.

\begin{lem} \label{estimativas} \mbox{}\\
\noindent (a) \, There exist $\epsilon_0,\sigma>0$ such that
$$
D_{\epsilon} \geq m_\lambda(c_0)+\sigma, \quad \forall \epsilon \in (0,\sigma).
$$
\noindent (b) \, $ \displaystyle \limsup_{\epsilon \to 0}\left\{\sup_{x \in X}
I_\epsilon(\phi_{\epsilon,x}) \right\} < 2m_\lambda(c_0)-\sigma$. \\
\noindent (c) \, There exist $\epsilon_0,R>0$ such that
$$
I_\epsilon(\phi_{\epsilon,x}) \leq \frac{1}{2}(m_\lambda(c_0)+D_\epsilon), \quad \forall \epsilon \in 
(0,\epsilon_0) \quad \mbox{and} \quad \forall x \in \partial B_R(0) \cap X.
$$ 
\end{lem}
\noindent {\bf Proof.} {\it  Proof of  (a),}  From definition of $D_{\epsilon} $, we know that
$$
D_{\epsilon}  \geq m_\lambda(c_0), \quad \forall \epsilon >0. 
$$
Supposing by contradiction that the lemma does not hold, there exists $\epsilon_n \to 0$ such that  
$$
D_{\epsilon_n} \to m_\lambda(c_0) \quad \mbox{as} \quad n \to +\infty.
$$
Hence, there exists $u_n \in {\mathcal N}_{\epsilon_n}$, with $\beta(u_n)\in Y$, satisfying 
$$
I_{\epsilon_n}(u_n) \to m_\lambda(c_0) \quad \mbox{as} \quad n \to +\infty.
$$
Thereby, by Lemma \ref{compacidade}, there exist $u_1 \in H^{s}(\mathbb{R}^{N}) \setminus \{0\}$ and $(z_n) 
\subset \mathbb{R}^N$  with $\displaystyle \liminf_{n \to +\infty}|\epsilon_n z_n|>0$ verifying
$$
u_n(\cdot+z_n) \to u_1 \quad \mbox{in} \quad  H^{s}(\mathbb{R}^{N}),
$$
that is,
$$
u_n=u_1(\cdot-z_n)+w_n \quad \mbox{with} \quad w_n \to 0 \quad \mbox{in} \quad H^{s}(\mathbb{R}^{N}).
$$
From definition of $\beta$, 
$$
\beta(u_1(\cdot-z_n))= \frac{\displaystyle  \int_{\mathbb{R}^{N}}\frac{\epsilon_n x+ \epsilon_n z_n}{|\epsilon_n x+
 \epsilon_n z_n|}|u_1|^{2}\,dx}{\displaystyle \int_{\mathbb{R}^{N}}|u_1|^{2}\,dx}.
$$
Repeating the same arguments explored in the proof of Lemma \ref{B1} (see also \cite{Benci}), we see that
$$
\beta(u_1(\cdot-z_n))=\frac{z_n}{|z_n|}+o_{n}(1),
$$
and so, 
$$
\beta(u_n)=\beta(u_1(\cdot-z_n))+o_n(1)=\frac{z_n}{|z_n|}+o_{n}(1).
$$
Since $\beta(u_n) \in Y$, we infer that $\frac{z_n}{|z_n|} \in Y_\lambda$ for $n$ large enough. Consequently,
 $z_n \in Y_\lambda$ for $n$ large enough, implying that
$$
\liminf_{n \to \infty}V(\epsilon_n z_n) > c_0.
$$
Making  $A=\displaystyle \liminf_{n \to \infty}V(\epsilon_n z_n)$, the last inequality together with Fatous's Lemma yields 
$$
m_\lambda(c_0)=\liminf_{n \to \infty}I_{\epsilon_n}(u_n) \geq J_{\lambda,A}(u_1) \geq  m_\lambda(A)>m_\lambda(c_0),
$$
which is an absurd, recalling that  $J'_A(u_1)u_1=0$ and $u_1 \not= 0$. \\

\noindent {\it Proof of (b).} \, Using condition $(V_4)$ , since  $u_0$ is a ground state solution
 associated with $I_{c_0}$, we deduce that there is $\epsilon_0>0$ such that 
$$
\sup_{x \in X}I_{\epsilon}(\phi_{\epsilon,x})\leq I_{c_0}(u_0)+\frac{3}{5}I_{c_0}(u_0)=m_\lambda(c_0)+ \frac{3}{5}m_\lambda(c_0)<2m_\lambda(c_0),  \quad \forall \epsilon \in (0, \epsilon_0).  
$$ 
\noindent {\it Proof of (c).} \, From $(V_1)$, given $\delta >0$, there are $R,\epsilon_0>0$ such that
$$
\sup\{I_\epsilon(\phi_{\epsilon,x})\,:\, x \in \partial B_{R}(0) \cap X \}\leq m_\lambda(c_0)+\delta \quad \forall \epsilon
 \in (0, \epsilon_0). 
$$ 
Fixing $\delta = \frac{\sigma}{4}$, where $\sigma$ was given in (a), we  have that
$$
\sup\{I_\epsilon(\phi_{\epsilon,x})\,:\, x \in \partial B_{R}(0) \cap X \}\leq \frac{1}{2}
\left(2m_\lambda(c_0)+\frac{\sigma}{2}\right)< \frac{1}{2}(m_\lambda(c_0)+D_\epsilon) \quad \forall \epsilon \in (0, \epsilon_0). 
$$

\fim

\vspace{0.5 cm}

Now, we are ready to show the minimax level. Define the map $\Phi_{\epsilon}:X \to H^{s}(\mathbb{R}^N)$ as  $\Phi_{\epsilon}=
\phi_{\epsilon,x}$. Denoting by $P$  the cone of nonnegative functions of $H^{s}(\mathbb{R}^{N})$, let us consider the set
$$
\Sigma=\{\Phi_\epsilon\,:\, x \in X, \, |x|\leq R \} \subset P,
$$
the class of functions
$$
\mathcal{H}=\left\{h \in C(P \cap \mathcal{N}_{\epsilon},P \cap \mathcal{N}_{\epsilon})\,:\,h(u)=u, \,\, \mbox{if} \,\, I_\epsilon(u) < \frac{1}{2}(m_\lambda(c_0)+D_\epsilon) \right\} 
$$
and finally the class of sets
$$
\Gamma=\{A \subset P \cap \mathcal{N}_\epsilon\,:\, A=h(\Sigma), \, h \in \mathcal{H} \}.
$$
\begin{lem} \label{intersecao}
If $A \in \Gamma$, then $A \cap \mathcal{B}_\epsilon \not= \emptyset$ for all $\epsilon>0$.
\end{lem}
\noindent {\bf Proof.} It is enough to show that for all $h \in  \mathcal{H}$, there is $x_* \in X$ with $|x_*| \leq R$ such that
$$
\beta(h(\Phi_\epsilon(x_*))) \in Y.
$$
For each $h \in \mathcal{H}$, we set the function $g:\mathbb{R}^{N} \to \mathbb{R}^{N}$ given by
$$
g(x)=\beta(h(\Phi_\epsilon(x))) \quad \forall x \in \mathbb{R}^{N},
$$
and the homotopy $\mathcal{F}:[0,1] \times X \to X$ as
$$
\mathcal{F}(t,x)=tP_X(g(x))+(1-t)x,
$$
where $P_X$ is the projection onto $X=\{(x,0)\,:\, x \in \mathbb{R}^N\}$. By using Corollary \ref{cor1}, fixed $R>0$ and $\epsilon>0$ small enough, we have that
$$
(\beta(\mathcal{F}(t,x)),x)>0, \quad \forall (t,x)  \in [0,1] \times ( \partial B_R \cap X).
$$
Using the homotopy invariance property of the Topological degree, we derive
$$
d(g,B_R \cap X,0)=1,
$$
implying that there exists $x_* \in B_R \cap X$ such that  $ \beta(h(\Phi_\epsilon(x_*)))=0$.
\fim

Now,  define the min-max value
$$
C_\epsilon=\inf_{A \in \Gamma}\sup_{u \in A}I_{\epsilon}(u).
$$
From Lemma \ref{intersecao},
\begin{equation} \label{NM1}
C_\epsilon \geq D_\epsilon \geq m_\lambda(c_0)+\sigma,
\end{equation}
for $\epsilon$ is small enough. On the other hand, 
$$
C_\epsilon \leq \sup_{x \in X}I_\epsilon(\phi_{\epsilon,x}), \quad \forall \epsilon >0.
$$ 
Then, by Lemma \ref{estimativas}(b), if $\epsilon$ is small enough
\begin{equation} \label{NM2}
C_\epsilon \leq \sup_{x \in X}I_\epsilon(\phi_{\epsilon,x})< 2m_\lambda(c_0)-\sigma.
\end{equation}
From (\ref{NM1}) and (\ref{NM2}), there is $\epsilon_0>0$ such that
\begin{equation} \label{NM3}
C_\epsilon \in (m_\lambda(c_0)+\sigma,2m_\lambda(c_0)-\sigma), \quad \forall \epsilon \in (0, \epsilon_0). 
\end{equation}

Now, we can use standard min-max arguments to conclude that $I_\epsilon$ has at least a critical point in $P \cap \mathcal{N}_{\epsilon}$ if $\epsilon$ is small enough.

\section{Final Comments}

The same approach used in the present paper can be used to prove the existence of solution for problems with subcritical growth like
$$
\epsilon^{2s} (-\Delta)^{s}{u}+V(z)u= f(u)\,\,\, \mbox{in} \,\,\, \mathbb{R}^{N}.
$$
To do this, it is enough to assume that $f$ verifies the following conditions:
\begin{enumerate}
\item[$(f_1)$] $\displaystyle \lim_{s\rightarrow
0}\dfrac{f(s)}{s}=0$.
\item[$(f_2)$] There is $\theta>2$ such that
$$
0<\theta F(s):=\theta\int_0^{s}f(t)dt\leq sf(s),\ \ \mbox{for all}\ \ s \in \mathbb{R} \setminus \{0\}. 
$$
\item[$(f_3)$] The function $s\rightarrow\dfrac{f(s)}{s}$ is
strictly increasing in $|s|>0$. 
\end{enumerate}

Related to function $V$ we assume $(V_1)-(V_4)$, however $(V_4)$ must be written of the following way  
$$ 
m(V(0)) \geq 2m(c_0) \quad \mbox{and} \quad c_1 \leq \left[1 + \frac{3}{5}\left(\frac{1}{2}-\frac{1}{\theta}\right)\right]c_0. \leqno{(V_4)} 
$$

\begin{acknowledgments}
This paper was completed  while the second author named was visiting the Department of Mathematics of the Rutgers University,
 whose hospitality he gratefully acknowledges. He would like to express his gratitude to Professor Haim Brezis  and Professor Yan Yan Li for invitation and friendship.

\end{acknowledgments}


\end{document}